\documentclass[runningheads,a4paper]{llncs}
\usepackage{amssymb}
\usepackage{graphicx}
\usepackage{url}
\newcommand{\keywords}[1]{\par\addvspace\baselineskip
\noindent\keywordname\enspace\ignorespaces#1}
\DeclareMathAlphabet{\mathpzc}{OT1}{pzc}{m}{it}
\newcommand{\mcA}{\mathcal{A}}

\newcommand{\mcD}{\mathcal{D}}

\newcommand{\mcL}{\mathcal{L}}

\newcommand{\mcN}{\mathcal{N}}
\newcommand{\mcQ}{\mathcal{Q}}
\newcommand{\ul}{\underline}
\newcommand{\ol}{\overline}

\newcommand{\diag}{\mbox{\textnormal{diag}}}
\newcommand{\tri}{\mbox{\textnormal{trig}}}
\newcommand{\nnz}{\mbox{\textnormal{nnz}}}

\usepackage{amssymb,amsmath,amsbsy}
\usepackage{enumitem}
\usepackage{graphicx}
\usepackage{algorithm2e}

\begin{document}
\mainmatter
\title{Efficient smoothers for all-at-once multigrid methods for Poisson and
Stokes control problems}

\titlerunning{Smoothers for Poisson and Stokes control problems}

%
%
\author{Stefan Takacs%
\thanks{The research was funded by the Austrian Science Fund (FWF): J3362-N25.}
}%
\authorrunning{Stefan Takacs}

\institute{Mathematical Institute,
University of Oxford,
United Kingdom\\
\path{stefan.takacs@numa.uni-linz.ac.at}\\
\url{http://www.numa.uni-linz.ac.at/~stefant/J3362/}}

%
%

\maketitle

\begin{abstract}
In the present paper we concentrate on an important
issue in constructing a good multigrid solver: the choice of an
efficient smoother. We will introduce all-at-once multigrid solvers
for optimal control problems which show robust convergence 
in the grid size and in the regularization parameter. We will refer
to recent publications that guarantee such a convergence behavior.
These publications do not pay much attention
to the construction of the smoother and suggest to use a normal
equation smoother. We will see that using
a Gauss Seidel like variant of this smoother, the overall
multigrid solver is speeded up by a factor of about two with
no additional work. The author will
give a proof which indicates that also the Gauss Seidel like variant 
of the smoother is covered by the convergence theory. Numerical
experiments suggest that the proposed method are competitive with 
Vanka type methods.
\end{abstract}

\section{Introduction}\label{sec:1}

In the present paper we discuss the construction of the all-at-once
multigrid solvers for two model problems.
The first model problem is a standard \emph{Poisson control problem}: 
Find a state $y \in H^1(\Omega)$ and a control $u\in L^2(\Omega)$
such that they minimize the cost functional
\begin{equation} \nonumber
   J(y,u) := \tfrac{1}{2} \|y-y_D\|_{L^2(\Omega)}^2 + \tfrac{\alpha}{2} \|u\|_{L^2(\Omega)}^2,
\end{equation}
subject to the elliptic boundary value problem (BVP)
\begin{equation}\nonumber
   -\Delta y + y = u \mbox{ in } \Omega
   \qquad \mbox{and} \qquad
     \tfrac{\partial y}{\partial n} = 0 \mbox{ on } \partial \Omega.
\end{equation}
The desired state $y_D$ and the regularization parameter $\alpha>0$ are assumed to be given. 
Here and in what follows, $\Omega\subseteq\mathbb{R}^2$ is a polygonal domain. We
want to solve the finite element discretization of this problem using a fast linear
solver which shows robust convergence behavior in the grid size and the
regularization parameter.
For solving this problem, we use the method of Lagrange multipliers, cf.~\cite{Lions:1971,Schoeberl:Simon:Zulehner:2010}.
We obtain a linear system in the state $y$, the control $u$ and the Lagrange multiplier $\lambda$. In
this linear system we eliminate the control as this has been done in~\cite{Schoeberl:Simon:Zulehner:2010,Takacs:Zulehner:2012}.
We discretize the resulting system using the Courant element and obtain a linear system:
\begin{equation}\label{eq:mp2}
     \underbrace{ \left(
	\begin{array}{cc}
	M_k & K_k \\
	K_k & -\alpha^{-1} M_k
        \end{array}
      \right)}_{\displaystyle\mcA_k:=}
      \underbrace{\left(
	\begin{array}{c}
	\ul{y}_k \\
	\ul{\lambda}_k
        \end{array}
      \right)}_{\displaystyle\ul{x}_k:=}
      =
      \underbrace{\left(
	\begin{array}{c}
	\ul{\mathpzc{f}}_k \\
	0
        \end{array}
      \right)}_{\displaystyle\ul{ f }_k:=}.
\end{equation}
Here, $M_k$ and $K_k$ are the standard
mass and stiffness matrices, respectively. The control can be
recovered using the following simple relation from the Lagrange multiplier: 
$\ul{u}_k = \alpha^{-1} \ul{\lambda}_k$, cf.~\cite{Schoeberl:Simon:Zulehner:2010}.
In \cite{Schoeberl:Simon:Zulehner:2010,Zulehner:2010} it was shown that
there are constants $\ul{C}>0$ and $\ol{C}>0$ (independent of the grid size
$h_k$ and the choice of $\alpha$) such that the stability estimate
\begin{equation}\label{eq:a1}
	\|\mathcal{Q}_k^{-1/2} \mcA_k \mathcal{Q}_k^{-1/2}\|\le \ol{C}\qquad \mbox{and}\qquad
	\|\mathcal{Q}_k^{1/2} \mcA_k^{-1} \mathcal{Q}_k^{1/2}\|\le \ul{C}^{-1}
\end{equation}
holds for the symmetric and positive definite matrix
\begin{equation*}
	\mathcal{Q}_k :=  \left(
		\begin{array}{cc}
			M_k + \alpha^{1/2}  K_k  \\
			  & \alpha^{-1} M_k + \alpha^{-1/2} K_k \\
		\end{array} 
	\right).
\end{equation*}

The second model problem is a standard \emph{Stokes control problem} (velocity tracking problem):
Find a velocity filed $v\in [H^1(\Omega)]^d$, a pressure distribution $p\in L^2(\Omega)$ and
a control $u\in [L^2(\Omega)]^d$ such that
\begin{equation}\nonumber
        J(v,p,u) = \tfrac12 \|v-v_D\|_{L^2(\Omega)}^2 + \tfrac{\alpha}{2} \|u\|_{L^2(\Omega)}^2
\end{equation}
is minimized subject to the Stokes equations
\begin{equation}\nonumber
        -\Delta v + \nabla p = u \mbox{ in } \Omega,  \qquad
        \nabla \cdot v = 0 \mbox{ in } \Omega, \qquad
        v = 0
        \mbox{ on } \partial \Omega. 
\end{equation}
The regularization parameter~$\alpha>0$ and the desired state (desired velocity field)
$v_D\in [L^2(\Omega)]^d$ are assumed to be given. To enforce uniqueness of the solution,
we additionally require $\int_{\Omega} p \mbox{ d}x=0$.

Similar as above, we can set up the optimality system and eliminate the control, 
cf.~\cite{Zulehner:2010,Takacs:2013a}. The
discretization can be done using the Taylor-Hood element. After these steps, we end up with
the following linear system:
\begin{equation}\label{eq:mp3}
      \underbrace{\left(
	\begin{array}{cccc}
	M_k & & K_k & D_k^T \\
	&0&D_k\\
	K_k &D_k^T& -\alpha^{-1} M_k\\
	D_k&&&0
        \end{array}
      \right)}_{\displaystyle\mcA_k:=}
      \underbrace{\left(
	\begin{array}{c}
	\ul{v}_k \\
	\ul{p}_k \\
	\ul{\lambda}_k \\
	\ul{\mu}_k \\
        \end{array}
      \right)}_{\displaystyle\ul{x}_k:=}
      =
      \underbrace{\left(
	\begin{array}{c}
	\ul{\mathpzc{f}}_k \\
	0\\0\\0
        \end{array}
      \right).}_{\displaystyle\ul{ f }_k:=}
\end{equation}
where $M_k$ and $K_k$ are standard mass and stiffness matrices and 
$D_k^T$ is the discretization of the gradient operator,
see, e.g.,~\cite{Zulehner:2010,Takacs:2013a}. Again, we are interested in
a fast solver which is robust in the regularization parameter and the grid size.
As in the previous example, the control $\ul{u}_k$ can by recovered from 
the Lagrange multiplier: $\ul{u}_k=\alpha^{-1}\ul{\lambda}_k$.
In~\cite{Zulehner:2010} it was shown that stability estimate~\eqref{eq:a1} is satisfied for 
\begin{equation*}
	\mathcal{Q}_k =  \mbox{block-diag}\left(
				W_k,\;
				\alpha D_kW_k^{-1}D_k^T,\;
				\alpha^{-1} W_k,\; 
				D_kW_k^{-1}D_k^T
		\right), 
\end{equation*}
where $W_k:=M_k + \alpha^{1/2}  K_k$.

\section{An all-at-once multigrid method}\label{sec:3}

The linear systems~\eqref{eq:mp2} and \eqref{eq:mp3} shall be solved
by a multigrid method, which reads as follows.
Starting from an initial approximation~$\ul{x}^{(0)}_k$,
one iterate of the multigrid method is given by the following two steps:
\begin{itemize}
        \item \emph{Smoothing procedure:} Compute
              \begin{equation} \nonumber
                   \ul{x}^{(0,m)}_k := \ul{x}^{(0,m-1)}_k + \hat{\mcA}_k^{-1}
                                    \left(\ul{ f}_k -\mcA_k\;\ul{x}^{(0,m-1)}_k\right)
                                    \qquad \mbox{for } m=1,\ldots,\nu
              \end{equation}
        with $\ul{x}^{(0,0)}_k=\ul{x}^{(0)}_k$. The choice of
        the smoother (or, in other words, of the  matrix
$\hat{\mcA}_k^{-1}$) will be discussed below. \vspace{.2cm}
        \item \emph{Coarse-grid correction:}
                \begin{itemize}
                     \item Compute the defect 
                        $\ul{ f}_k -\mcA_k\;\ul{x}^{(0,\nu)}_k$
                        and restrict it to grid level $k-1$ using
                        an restriction matrix $I_k^{k-1}$:\;
                       $
                              \ul{r}_{k-1}^{(1)} := I_k^{k-1} \left(\ul{ f}_k -\mcA_k
                              \;\ul{x}^{(0,\nu)}_k\right).
                       $
                     \item Solve the following coarse-grid problem approximatively:
                        \begin{equation}\label{eq:coarse:grid:problem}
                            \mcA_{k-1} \,\ul{p}_{k-1}^{(1)} =\ul{r}_{k-1}^{(1)}
                        \end{equation}
                     \item Prolongate $\ul{p}_{k-1}^{(1)}$  to the
                          grid level $k$ using an prolongation 
                          matrix $I^k_{k-1}$ and add
                          the result to the previous iterate:
                          $
                               \ul{x}_{k}^{(1)} := \ul{x}^{(0,\nu)}_k +
                                I_{k-1}^k \, \ul{p}_{k-1}^{(1)}.
                          $
                \end{itemize}
\end{itemize}
As we have assumed to have nested spaces, the intergrid-transfer
matrices can be chosen in
a canonical way:
$I_{k-1}^k$ is the canonical embedding and the restriction $I_k^{k-1}$
is its (properly scaled) transpose.
If the problem~\eqref{eq:coarse:grid:problem} is solved exactly,
we obtain the two-grid method.
In practice, the problem~\eqref{eq:coarse:grid:problem} is
approximatively solved by applying one step (V-cycle)
or two steps (W-cycle) of the multigrid method, recursively. Only
the coarsest grid level,~\eqref{eq:coarse:grid:problem} is 
solved exactly.

The only part  of the multigrid algorithm that has not been specified yet,
is the smoother. For the choice of the smoother, we make use of the convergence
theory. We develop a convergence theory based on Hackbusch's
splitting of the analysis into smoothing property and approximation
property:
\begin{itemize}
        \item \emph{Smoothing property:}
                \begin{equation} \label{eq:smp}
                         \sup_{\ul{\tilde{x}}_k\in X_k} 
                         \frac{\left(\mcA_k(\ul{x}_k^{(0,\nu)}-\ul{x}_k^*),
                                   \ul{\tilde{x}}_k\right)_{\ell^2}}{\| \ul{\tilde{x}}_k\|_{\mcL_k}}
                         \le \eta(\nu) \| \ul{x}_k^{(0)}-\ul{x}_k^*\|_{\mcL_k}
                \end{equation}
                should hold for some function $\eta(\nu)$
                with $\lim_{\nu\rightarrow\infty}\eta(\nu)= 0$. Here and in what follows,
					$\ul{x}_k^* := \mcA_k^{-1} \ul{ f }_k$ is the exact solution, $\|\cdot\|_{\mcL_k}:=
					(\cdot,\cdot)_{\mcL_k}^{1/2} := (\mcL_k \cdot,\cdot)_{\ell^2}^{1/2}$ for some
					symmetric positive definite matrix $\mcL_k$ and 
					$(\cdot,\cdot)_{\ell^2}$ is the standard Euclidean scalar product.
        \item \emph{Approximation property:}
                \begin{equation} \nonumber
                        \| \ul{x}_k^{(1)}-\ul{x}_k^*\|_{\mcL_k}\le
                        C_A \sup_{\ul{\tilde{x}}_k\in X_k} 
                        \frac{\left(\mcA_k(\ul{x}_k^{(0,\nu)}-\ul{x}_k^*),
                                   \ul{\tilde{x}}_k\right)_{\ell^2}}{\| \ul{\tilde{x}}_k\|_{\mcL_k}}
                \end{equation}        
                should hold for some constant $C_A>0$.
\end{itemize}
It is easy to see that, if we combine both conditions, we see 
that the two-grid method converges in the norm $\|\cdot\|_{\mcL_k}$ for $\nu$ large enough.
The convergence of the W-cycle multigrid method can be shown under
mild assumptions, see e.g.~\cite{Hackbusch:1985}.

For the smoothing analysis, it is convenient to rewrite the smoothing property in pure matrix notation:
\eqref{eq:smp} is equivalent to
\begin{equation} \label{eq:smp2}
	\|\mcL_k^{-1/2}\mcA_k(I-\hat{\mcA}_k^{-1}\mcA_k)^{\nu}\mcL_k^{-1/2}\| \le \eta(\nu).
\end{equation}
For the Poisson control problem, it was shown in \cite{Schoeberl:Simon:Zulehner:2010},
that the approximation property is satisfied for the following choice of the matrix $\mcL_k$ (note that this matrix
represents the norm $\|\cdot\|_{X^-}$ used in the mentioned paper)
\begin{equation*}
	\mcL_k =  \left(
		\begin{array}{cc}
			\diag(M_k + \alpha^{1/2}  K_k)  \\
			& \diag(\alpha^{-1} M_k + \alpha^{-1/2} K_k) \\
		\end{array} 
		\right),
\end{equation*}
i.e., $\mcL_k = \diag(\mcQ_k)$. Here and in what follows, $\diag(M)$ is the diagonal matrix containing the diagonal of a matrix $M$.
For the Stokes control problem it was shown in \cite{Takacs:2013a}, that the approximation
property is satisfied for the following choice of $\mcL_k$: 
\begin{equation*}
	\mcL _k = \left(
		\begin{array}{cccc}
			\hat{W}_k \\ & \hat{P}_k \\ && \alpha^{-1} \hat{W}_k \\ &&&\alpha^{-1} \hat{P}_k 
		\end{array} 
	\right),
\end{equation*}
where $\hat{W}_k := \diag(M_k+\alpha^{1/2} K_k)$ and $\hat{P}_k := \alpha \;\diag( D_k \hat{W}_k^{-1} D_k^T ).$

Still, we have not specified the choice of the smoother, which now can be done using
the convergence theory. We have seen for which choices of $\mcL_k$ the approximation
property is satisfied. We are interested in a smoother such that the smoothing
property is satisfied for the same choice of $\mcL_k$.

In~\cite{Takacs:Zulehner:2012,Takacs:2013a} a
\emph{normal equation smoother} was proposed. This approach
is applicable to a quite general class of problems, cf.~\cite{Brenner:1996}
and others. In our notation, the normal equation smoother reads as follows:
\begin{equation}\nonumber
        \ul{x}^{(0,m)}_k := \ul{x}^{(0,m-1)}_k + \tau
                         \underbrace{\mcL_k^{-1} \mcA_k^T \mcL_k^{-1}}_{\displaystyle
\hat{\mcA}_k^{-1}:=}
                        \left(\ul{ f}_k -\mcA_k \;\ul{x}^{(0,m-1)}_k\right)
                \quad \mbox{for } m=1,\ldots,\nu.
\end{equation}
Here, a fixed~$\tau>0$ has to be chosen such that the
spectral radius~$\rho(\tau \hat{\mcA}_k^{-1}\mcA_k)$ is bounded away
from~$2$ on all grid levels~$k$ and for all choices of the parameters.
It was shown that it is possible to find such an uniform $\tau$
for the Poisson control problem, e.g., in \cite{Takacs:Zulehner:2012} and
for the Stokes control problem, e.g., in \cite{Takacs:2013a}. For
the normal equation smoother, the smoothing property can be shown using
a simple eigenvalue analysis, cf.~\cite{Brenner:1996}.
Numerical experiments show that the normal equation smoother works rather well for the
mentioned model problems. However, there are smoothers such that the overall
multigrid method converges much faster. Note
that the normal equation smoother is basically a Richardson iteration scheme, applied
to the normal equation. It is well-known for elliptic problems that Gauss Seidel
iteration schemes are typically much better smoothers than Richardson iteration schemes.
In the context of saddle point problems, the idea of Gauss Seidel
smoothers has been applied, e.g., in the context of collective smoothers, see below. However, in the context
of normal equation smoothers the idea of Gauss Seidel smoothers has not gained much
attention. The setup of such an approach is straight forward: In compact notation 
such an approach, which we call \emph{least squares Gauss Seidel} (LSGS) approach, 
reads as follows:
\begin{equation}\nonumber
        \ul{x}^{(0,m)}_k := \ul{x}^{(0,m-1)}_k + 
                        \underbrace{ \tri(\mcN_k)^{-1} \mcA_k^T \mcL_k^{-1}}_{\displaystyle\hat{\mcA}_k:=}
                        \left(\ul{ f}_k -\mcA_k \;\ul{x}^{(0,m-1)}_k\right)
                \;\; \mbox{for } m=1,\ldots,\nu,
\end{equation}
where $\mcN _k:=\mcA_k^T\mcL_k^{-1} \mcA_k$ and $\tri(M)$ is a matrix whose coefficients
coincide with the coefficients of $M$ on the diagonal and the left-lower triangular part and
vanish elsewhere.
The author provides a possible realization of that approach as Algorithm~\ref{alg2} to
convince the reader that the computational complexity of the LSGS approach is equal
to the computational complexity of the normal equation smoother, where a possible
realization is given as Algorithm~\ref{alg1}.

We will see below that the LSGS approach works very well in the numerical experiments. However,
there is no proof of the smoothing property known to the author. This is 
due to the fact that the matrix $\hat{\mcA}_k$ is not symmetric. 
One possibility to overcome this difficulty is to consider the symmetric version (symmetric
least squares Gauss Seidel approach, sLSGS approach). 
This is analogous to the case of elliptic problems: For elliptic problems the smoothing
property for the symmetric Gauss Seidel iteration can be shown for general cases but
for the standard Gauss Seidel iteration the analysis is restricted to special cases, 
cf. Section~6.2.4 in~\cite{Hackbusch:1985}.

\begin{algorithm}[t]
	\textbf{Given:} Iterate $(\texttt{x}_i)_{i=1}^{N}=\ul{x}^{(0,m-1)}$ and corresp. residual $(\texttt{r}_i)_{i=1}^{N}=\ul{f} - \mcA \ul{x}^{(0,m-1)}$;\\
	\textbf{Result:} Iterate $(\texttt{x}_i)_{i=1}^{N}=\ul{x}^{(0,m)}$ and corresp. residual $(\texttt{r}_i)_{i=1}^{N}=\ul{f} - \mcA \ul{x}^{(0,m)}$;\\
	\For{$i=1,\ldots,N$}{
		$\texttt{q} := 0$;\\
		\mbox{\textnormal{\textbf{for all $j$ such that $\mcA_{i,j}\not=0$}}} \mbox{\textnormal{\textbf{do}}}
			$\texttt{q} := \texttt{q} + \mcA_{i,j} / \mcL_{j,j} * \texttt{r}_j$;\\
		$\texttt{p}_i := \tau * \texttt{q} / \mcL_{i,i}$;
	}
	\For{$i=1,\ldots,N$}{
		$\texttt{x}_i := \texttt{x}_i + \texttt{p}_i$;\\
		\mbox{\textnormal{\textbf{for all $j$ such that $\mcA_{j,i}\not=0$}}} \mbox{\textnormal{\textbf{do}}}
			$\texttt{r}_j := \texttt{r}_j - \mcA_{j,i} * \texttt{p}_i $;
	}
	\caption{Normal equation iteration scheme}\label{alg1}
\end{algorithm}
\begin{algorithm}[t]
	\textbf{Given:} Iterate $(\texttt{x}_i)_{i=1}^{N}=\ul{x}^{(0,m-1)}$ and corresp. residual $(\texttt{r}_i)_{i=1}^{N}=\ul{f} - \mcA \ul{x}^{(0,m-1)}$;\\
	\textbf{Result:} Iterate $(\texttt{x}_i)_{i=1}^{N}=\ul{x}^{(0,m)}$ and corresp. residual $(\texttt{r}_i)_{i=1}^{N}=\ul{f} - \mcA \ul{x}^{(0,m)}$;\\
	\textbf{Prepare once:} $\mcN _{i,i}:=\sum_{j=1}^{N} \mcA_{i,j}^2 / \mcL_{j,j}$ for all $i=1,\ldots,N$;\\
	\For{$i=1,\ldots,N$}{
		$\texttt{q} := 0$;\\
		\mbox{\textnormal{\textbf{for all $j$ such that $\mcA_{i,j}\not=0$}}} \mbox{\textnormal{\textbf{do}}}
			$\texttt{q} := \texttt{q} + \mcA_{i,j} / \mcL_{j,j} * \texttt{r}_j$;
		
		$\texttt{p} := \texttt{q} / \mcN _{i,i}$;\\
		$\texttt{x}_i := \texttt{x}_i + \texttt{p}$;\\
		\mbox{\textnormal{\textbf{for all $j$ such that $\mcA_{j,i}\not=0$}}} \mbox{\textnormal{\textbf{do}}}
			$\texttt{r}_j := \texttt{r}_j - \mcA_{j,i} * \texttt{p} $;
		
	}
	\caption{LSGS iteration scheme}\label{alg2}
\end{algorithm}

One step of the sLSGS iteration consists of one step of the LSGS iteration, followed by one step
of the LSGS iteration with reversed order of the variables. (So the computational complexity of
one step of the sLSGS iteration is equal to the computational complexity of two steps of the standard
LSGS iteration.) One step of the sLSGS iteration reads as follows in compact notation:
\begin{align}\nonumber
        &\ul{x}^{(0,m)}_k := \ul{x}^{(0,m-1)}_k + 
                        \hat{\mcN}_k^{-1} \mcA_k^T \mcL_k^{-1}
                        \left(\ul{ f}_k -\mcA_k \;\ul{x}^{(0,m-1)}_k\right)
			\qquad \mbox{for } m=1,\ldots,\nu,\\
		&\mbox{where $\hat{\mcN}_k:=\tri(\mcN_k)\; \diag(\mcN_k)^{-1}\; \tri(\mcN_k)^T$.}\label{eqmcNhat}
\end{align}
For our needs, the following convergence lemma is sufficient.
\begin{lemma}\label{lem}
	Assume that $\mcA_k$ is sparse, \eqref{eq:a1} is satisfied and let $\mcL_k$ be a positive definite
	diagonal matrix such that
	\begin{equation}\label{eq:lem}
		\|\mcQ_k^{1/2}\ul{x}_k\|\le \|\mcL_k^{1/2}\ul{x}_k\|\quad\mbox{for all } \ul{x}_k.
	\end{equation}
	Then the sLSGS approach satisfies the smoothing property~\eqref{eq:smp2}, i.e.,
	\begin{equation*}
		\|\mcL_k^{-1/2} \mcA_k (I-\hat{\mcN}_k^{-1} \mcN_k)^{\nu} \mcL_k^{-1/2}\| \le \frac{2^{-1/2}\;\ol{C}\;\nnz(\mcA_k)^{5/2} }{\sqrt{\nu}},
	\end{equation*}
	where $\nnz(M)$ is the maximum number of non-zero entries per row of~$M$.
\end{lemma}
Note that \eqref{eq:lem} is a standard inverse inequality, which is satisfied for both model problems, cf.~\cite{Schoeberl:Simon:Zulehner:2010,Takacs:Zulehner:2012,Takacs:2013a}.
Note moreover that this assumption also has to be satisfied to show the smoothing property for the normal equation smoother, cf.~\cite{Takacs:Zulehner:2012,Takacs:2013a}.

{\em Proof of Lemma~\ref{lem}.}
	The combination of \eqref{eq:a1} and \eqref{eq:lem} yields $\|\mcL_k^{-1/2}\mcA_k \mcL_k^{-1/2}\|\le \ol{C}$.
	Prop.~6.2.27 in~\cite{Hackbusch:1985} states that for any symmetric positive definite matrix $\mcN_k$
	\begin{equation}\label{eq:hb}
		\|\hat{\mcN}_k^{-1/2} \mcN_k (I-\hat{\mcN}_k^{-1} \mcN_k)^{\nu} \hat{\mcN}_k^{-1/2}\| \le \nu^{-1}
	\end{equation}
	holds, where $\hat{\mcN}_k$ is as in~\eqref{eqmcNhat}. Using $\mcD_k:= \mbox{diag}(\mcN_k)$, we obtain
	\begin{align*}
		&\|\mcL_k^{-1/2} \hat{\mcN}_k^{1/2}\|^2  = \rho( \mcL_k^{-1/2}\hat{\mcN}_k\mcL_k^{-1/2}) 
		\le \|\mcL_k^{-1/2}\tri(\mcN_k) \mcD_k^{-1/2}\|^2\\
		&\quad 
		 \le\|\mcL_k^{-1/2}\mcD_k^{1/2}\|^2 \|\mcD_k^{-1/2}\tri(\mcN_k) \mcL_k^{-1/2}\|^2
	\end{align*}
	Let $\mcA_k=(\mcA_{i,j})_{i,j=1}^N$, $\mcN_k=(\mcN_{i,j})_{i,j=1}^N$,
	 $\mcL_k=(\mcL_{i,j})_{i,j=1}^N$ and $\psi(i):=\{j\in\mathbb{N}:\mcN_{i,j}\not=0\}$. We obtain using 
	 Gerschgorin's theorem, the fact that the infinity norm is monotone in the matrix entries,
	and using the symmetry of $\mcN_k$ and $\mcA_k$ and Cauchy-Schwarz inequality:
	\begin{align}
		&\|\mcD_k^{-1/2}\tri(\mcN_k) \mcD_k^{-1/2}\| \nonumber\\
		& \le\|\mcD_k^{-1/2}\tri(\mcN_k) \mcD_k^{-1/2}\|_{\infty}^{1/2}
		\|\mcD_k^{-1/2}\tri(\mcN_k)^T \mcD_k^{-1/2}\|_{\infty}^{1/2} 
		 	 \le \|\mcD_k^{-1/2}\mcN_k\mcD_k^{-1/2}\|_{\infty}  \nonumber\\
			&  = \max_{i=1,\ldots, N}  \sum_{k\in\psi(i)}
					\left(\sum_{n=1}^N \frac{\mcA_{i,n}^2}{\mcL_{n,n}} \right)^{-1/2}
					\left(\sum_{j=1}^N \frac{\mcA_{i,j} \mcA_{j,k}}{\mcL_{j,j}}\right)
					\left(\sum_{n=1}^N \frac{\mcA_{k,n}^2}{\mcL_{n,n}} \right)^{-1/2}\nonumber\\
			&  \le \max_{i=1,\ldots, N}  \sum_{k\in\psi(i)} 1 = \nnz(\mcN_k) \le \nnz(\mcA_k)^2.\label{eq:xx1}
	\end{align}
	Further, we obtain
	\begin{align}
		&\|\mcL_k^{-1/2}\mcD_k^{1/2}\|^2 =\|\mcL_k^{-1/2}\mcD_k^{1/2}\|_{\infty}^2 = \|\mcL_k^{-1/2}\mcD_k\mcL_k^{-1/2}\|_{\infty}
		 = \max_{i=1,\ldots,N} \sum_{j=1}^N \frac{\mcA_{i,j}^2}{\mcL_{i,i}\mcL_{j,j}} \nonumber\\
		&\quad\le \nnz(\mcA_k) \max_{i,j=1,\ldots,N}\frac{\mcA_{i,j}^2}{\mcL_{i,i}\mcL_{j,j}} 
		= \nnz(\mcA_k) \|\mcL^{-1/2}\mcA\mcL^{-1/2}\|^2 
		\le \nnz(\mcA_k) \;\ol{C}^2.\label{eq:xx2}
	\end{align}
	By combining \eqref{eq:hb}, \eqref{eq:xx1} and \eqref{eq:xx2}, we obtain
	\begin{align*}
		&\|\mcL_k^{-1/2} \mcA_k (I-\hat{\mcN}_k^{-1} \mcN_k)^{\nu} \mcL_k^{-1/2}\|^2 \\
		& \quad \le \|\mcL_k^{-1/2}(I- \mcN_k\hat{\mcN}_k^{-1})^{\nu}\mcA_k\mcL_k^{-1} \mcA_k (I-\hat{\mcN}_k^{-1} \mcN_k)^{\nu} \mcL_k^{-1/2}\|\\
		& \quad = \|\mcL_k^{-1/2}\mcN_k(I- \hat{\mcN}_k^{-1}\mcN_k)^{2\nu} \mcL_k^{-1/2}\| \le \frac{\ol{C}^2 \nnz(\mcA_k)^5}{2\nu},
	\end{align*}
	which finishes the proof.
\qed

We went to compare the numerical behavior of the LSGS approach with the behavior of
a standard smoother. One class of standard smoothers for saddle point problems is the class 
of \emph{Vanka type smoothers}, which has been originally introduced
for Stokes problems, cf.~\cite{Vanka:1986}. Such smoothers have also gained interest
for optimal control problems, see, e.g.,~\cite{Trottenberg:2001,Borzi:Kunisch:Kwak:2003,Takacs:Zulehner:2011}.

The idea of Vanka type smoothers is to compute updates in subspaces directly for
the whole saddle point problem and to combine these updates is an additive or a
multiplicative way to compute the next update. Here, the variables are not grouped 
based on the block-structure of $\mcA_k$, but the
grouping is done of based on the location of the corresponding degrees of freedom in the
domain $\Omega$. The easiest of such ideas for the Poisson 
control problems is to do the grouping point-wise, which leads to the idea of \emph{point smoothing}.
Here, we group for each node $\delta_i$ of the discretization (each degree of freedom of the Courant element)
the value $y_i$ of the state and the value $\lambda_i$ of the Lagrange multiplier and
compute an update in the corresponding subspace. The multiplicative variant of such a smoother is
a \emph{collective Gauss Seidel} (CGS) smoother:
\begin{align*}
        \ul{x}^{(0,m,i)}_k &:= \ul{x}^{(0,m,i-1)}_k + 
                         \mathcal{P}_k^{(i)} \left(\left.\mathcal{P}_k^{(i)}\right.^T \mcA_k\mathcal{P}_k^{(i)}\right)^{-1} \left.\mathcal{P}_k^{(i)}\right.^T
                        \left(\ul{f}_k -\mcA_k \;\ul{x}^{(0,m,i-1)}_k\right),
\end{align*}
where $\ul{x}^{(0,m,0)}_k:=\ul{x}^{(0,m-1)}_k$  and $\ul{x}^{(0,m)}_k :=\ul{x}^{(0,m,N_k)}_k$.
For each $i=1,\ldots, N_k$, the matrix
$\mathcal{P}_k^{(i)}\in \mathbb{R}^{2 N_k\times 2}$ takes the value $1$ on the positions $(i,1)$ and
$(i+N_k,2)$ and the value $0$ elsewhere. For the Poisson control problem, we obtain
\begin{equation*}
	\left.\mathcal{P}_k^{(i)}\right.^T \mcA_k\mathcal{P}_k^{(i)} = \left(
			\begin{array}{cc}
					 M_{i,i} &  K_{i,i} \\ K_{i,i} & - \alpha^{-1} M_{i,i}
			\end{array}
		\right),
\end{equation*}
where $M_{i,i}$ and $K_{i,i}$ are the entries of the matrices
$M_k$ and $K_k$.

\begin{figure}
	\begin{center}
	\vspace{-.2cm}
	\includegraphics[scale=.3]{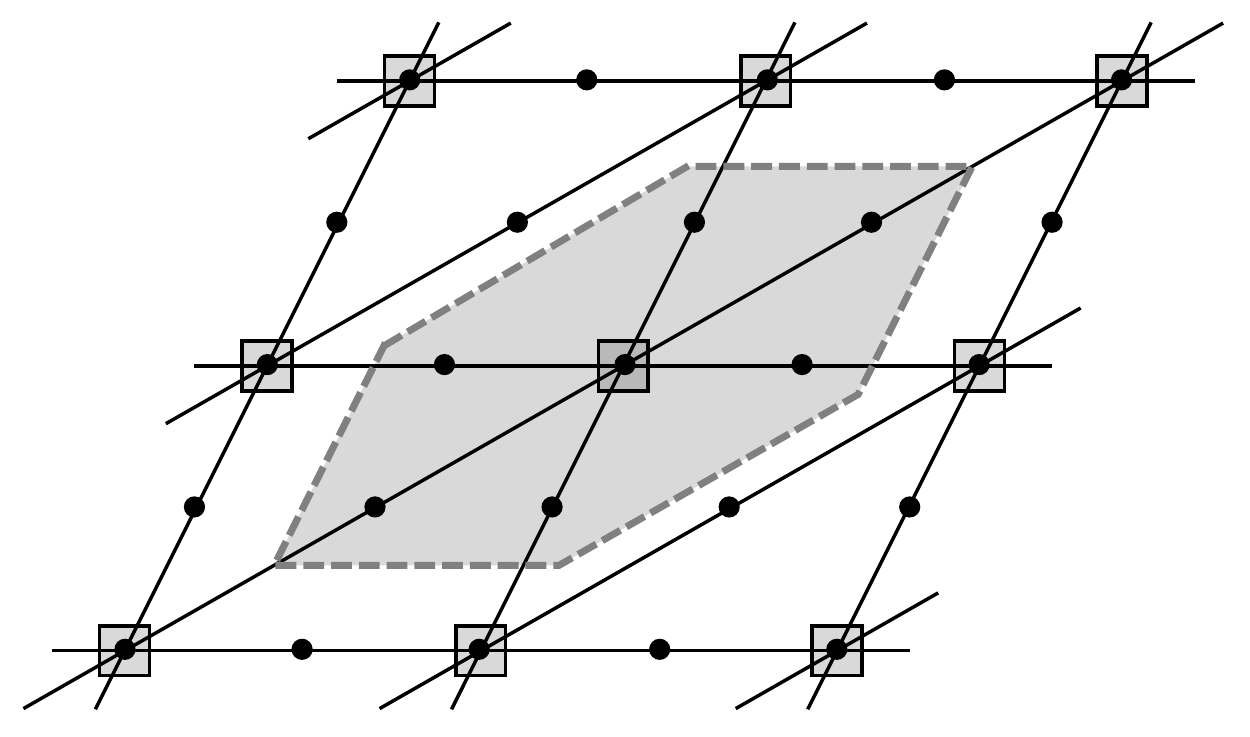}
	\caption{Patches for the Vanka-type smoother applied to a Taylor Hood discretization.
	The dots are the degrees of freedom of $v$ and $\lambda$,
	the rectangles are the degrees of freedom of $p$ and $\mu$}\label{fig1}
	\vspace{-.6cm}
	\end{center}
\end{figure}

For the Stokes control problem, it is not reasonable to use exactly the same approach.
This is basically due to the fact that the degrees of freedom for $v$ and $\lambda$ are not
located on the same positions as the degrees of freedom for $p$ and $\mu$. However,
we can introduce an approach based on patches: so, for each vertex of the triangulation,
we consider subspaces that consist of the degrees of 
freedoms located on the vertex itself and the degrees of freedom located on all edges
which have one end at the chosen vertex, cf. Fig.~\ref{fig1}. 
Note that here the subspaces are much larger
than the subspaces chosen in the case of the CGS approach for the Poisson control problem (which was just $2$).
This increases the computational cost of applying the method significantly.
For Vanka type smoothers there are only a few convergence results known,
cf.~\cite{Borzi:Kunisch:Kwak:2003} for a Fourier Analysis and an analysis based compactness argument
and \cite{Takacs:Zulehner:2011} for a proof based on Hackbusch's splitting of the analysis
into smoothing property and approximation property which shows the convergence in case of a collective
Richardson smoother.

\section{Numerical results}\label{sec:5}

In this section we give numerical results to illustrate quantitatively
the convergence behavior of the proposed methods.
The number of iterations was measured as
follows: We start with a random initial guess and iterate until
the relative error in the norm ${\|\cdot\|_{\mcL_k}}$ was reduced by
a factor of $10^{-6}$. Without loss of generality, the right-hand side
was chosen to be $0$.
For both model problems, the normal equation smoother,
the LSGS smoother, the sLSGS smoother and a Vanka type smoother
 have been applied. For the smoothers
$2$ pre- and $2$ post-smoothing steps have been applied. Only for the
sLSGS smoother, just $1$ pre- and $1$ 
post-smoothing step has been applied. This is due to the fact that one
step of the symmetric version is basically the same computational cost
as two steps of the standard version.
The normal equation smoother was damped with $\tau=0.4$ for the Poisson
control problem and $\tau=0.35$ for the Stokes control problem,
cf.~\cite{Takacs:Zulehner:2012,Takacs:2013a}.
For the Gauss Seidel-like approaches, damping was not used.

In Table~\ref{tab:1}, we give the results for the standard 
Poisson control problem. Here, we see that all smoothers lead to
convergence rates that are well bounded for a wide range of $h_k$ and $\alpha$. Compared
to the normal equation smoother, the LSGS smoother
leads to a speedup be a factor of about two without any additional work.
The symmetric version (sLSGS) is a bit slower than the LSGS method. 
For the first model problem, the (popular) CGS method is significantly faster. 
However, for this method no convergence theory is known. 
\begin{table}%
\begin{center}
        \begin{tabular}{p{1.cm}p{.6cm}p{.6cm}p{1.3cm}p{.6cm}p{.6cm}p{1.3cm}p{.6cm}p{.6cm}p{1.3cm}p{.6cm}p{.6cm}p{.8cm}}
        \hline\noalign{\smallskip}
               & \multicolumn{3}{l}{Normal equation} &
		 \multicolumn{3}{l}{LSGS} &
		 \multicolumn{3}{l}{sLSGS} &
                 \multicolumn{3}{l}{CGS} \\
        \noalign{\smallskip}\hline\noalign{\smallskip}
		$\alpha=$&$10^0$&$10^{-6}$&$10^{-12}$&$10^0$&$10^{-6}$&$10^{-12}$&$10^0$&$10^{-6}$&$10^{-12}$&$10^0$&$10^{-6}$&$10^{-12}$\\
        \noalign{\smallskip}\hline\noalign{\smallskip}
		$k=5$&26&31&28&11& 9& 7&14&12&14& 5& 5& 3\\
		$k=6$&27&28&29&11&11& 7&14&14&13& 5& 5& 3\\
		$k=7$&27&28&31&11&11& 6&14&14&12& 5& 5& 3\\
		$k=8$&27&27&25&11&11& 3&14&14& 7& 5& 5& 4\\
        \noalign{\smallskip}\hline\noalign{\smallskip}
        \end{tabular}
        \caption{Number of iterations for the \emph{Poisson control model problem}}
        \label{tab:1}      
	\vspace{-.6cm}
\end{center}
\end{table}%

In Table~\ref{tab:3}, we give the convergence results for the Stokes
control problem. Also here we observe that the LSGS and the sLSGS approach
lead to a speedup of a factor of about two compared to the normal
equation smoother. Here, the Vanka type smoother shows slightly
smaller iteration numbers than the LSGS approach. In terms of computational
costs, the LSGS smoother seems to be much better than the patch-based
Vanka type smoother because there relatively large subproblems have to be solved
to compute the updates.
This is different the case of the CGS smoother, where the
subproblems are just $2$-by-$2$ linear systems.
Numerical experiments have shown that the undamped version of the patch-based Vanka type method does
not lead to a convergent multigrid method. So, this smoother was damped with
$\tau=0.4$. Due to lack of convergence theory, the author cannot explain
why this approach -- although it is a multiplicative approach -- needs damping.
\begin{table}%
\begin{center}
	\vspace{-.2cm}
        \begin{tabular}{p{1.cm}p{.6cm}p{.6cm}p{1.3cm}p{.6cm}p{.6cm}p{1.3cm}p{.6cm}p{.6cm}p{1.3cm}p{.6cm}p{.6cm}p{.8cm}}
        \hline\noalign{\smallskip}
               & \multicolumn{3}{l}{Normal equation} &
		 \multicolumn{3}{l}{LSGS} &
		 \multicolumn{3}{l}{sLSGS} &
                 \multicolumn{3}{l}{Vanka type} \\
        \noalign{\smallskip}\hline\noalign{\smallskip}
		$\alpha=$&$10^0$&$10^{-6}$&$10^{-12}$&$10^0$&$10^{-6}$&$10^{-12}$&$10^0$&$10^{-6}$&$10^{-12}$&$10^0$&$10^{-6}$&$10^{-12}$\\
        \noalign{\smallskip}\hline\noalign{\smallskip}
		$k=4$&31&31&60&13&12&14&17&16&22&11&10& 7\\
		$k=5$&32&30&55&14&13&12&18&16&19&11&10& 7\\
		$k=6$&32&31&44&14&13& 9&18&17&12&11&11& 7\\
		$k=7$&32&31&37&14&14& 6&18&17& 9&11&11& 9\\
        \noalign{\smallskip}\hline\noalign{\smallskip}
        \end{tabular}
        \caption{Number of iterations for the \emph{Stokes control model problem}}
        \label{tab:3}  
	\vspace{-.6cm}     
\end{center}
\end{table}%

For completeness, the author wants to mention that for cases, where a (closed form of a) matrix $\mcQ_k$ 
satisfying~\eqref{eq:a1} robustly is not known, the normal equation smoother
does not show as good results as methods where such an information is not needed,
like Vanka type methods. This was discussed in~\cite{Takacs:Zulehner:2011} for a
boundary control problem, but it is also true for the linearization of optimal control
problems with inequality constraints as discussed in~\cite{Herzog:Sachs:2010}
and others. The same is true for the Gauss Seidel like variants of the normal equation smoother.

Concluding, we have observed that accelerating the idea of normal
equation smoothing with a Gauss Seidel approach, leads to a speedup
of a factor of about two without any further work. The fact that convergence
theory is known for the sLSGS approach, helps also for the numerical practice
(unlike the case of Vanka type smoothers).

\bibliographystyle{amsplain}
\bibliography{literature}
\vspace{-3cm}
\parbox{12cm}{ \vspace{3cm}
\begin{center}
This paper has been published in\\
C. P\"otzsche, C. Heuberger, B. Kaltenbacher and F. Rendl: \\System Modeling and Optimization. Springer, 2014.\\[1em]
The original publication is available at www.springerlink.com:\\
\url{http://link.springer.com/chapter/10.1007/978-3-662-45504-3_33}\end{center} \vspace{-3cm}}

\end{document}